\documentclass[11pt, preprint]{imsart}
\usepackage{xr}
\usepackage{amsmath,amssymb}
\usepackage{parskip}
\usepackage{marvosym}
\usepackage{mathtools,bm}
\usepackage[colorlinks]{hyperref}
\hypersetup{
	colorlinks,%
	citecolor=blue,%
	filecolor=black,%
	linkcolor=blue,%
	urlcolor=blue
}
\usepackage{enumerate}
\usepackage{multirow}
\RequirePackage[OT1]{fontenc}

\usepackage{amsthm}
\usepackage{amsmath}
\usepackage{amssymb}
\usepackage{thmtools} 
\usepackage{textcomp}

\usepackage{graphicx}
\usepackage{verbatim}
\usepackage{array, float}
\usepackage{fontenc}
\usepackage[toc,page]{appendix}
\usepackage[nottoc]{tocbibind}
\usepackage[english]{babel}
\usepackage{marvosym}

\usepackage[pagewise]{lineno}

\usepackage{mathtools}
\allowdisplaybreaks
\usepackage{bbm}
\usepackage[noabbrev,capitalize]{cleveref}

\newtheoremstyle{exampstyle}
{8pt} 
{8pt} 
{\it} 
{} 
{\bfseries} 
{.} 
{.5em} 
{} 

\theoremstyle{exampstyle}
\newtheorem{theorem}{Theorem}

\newtheorem{remark}{Remark}

\newtheorem{prop}{Proposition}

\newtheorem{defn}{Definition}

\numberwithin{equation}{section}
\numberwithin{example}{section}
\numberwithin{theorem}{section}
\numberwithin{lemma}{section}
\numberwithin{corollary}{section}
\numberwithin{prop}{section}
\numberwithin{defn}{section}
\numberwithin{remark}{section}

\usepackage{tikz}
\tikzset{
	treenode/.style = {shape=rectangle, rounded corners,
		draw, align=center,
		top color=white, bottom color=blue!20},
	root/.style     = {treenode, font=\Large, bottom color=yellow},
	env/.style      = {treenode, font=\ttfamily\normalsize},
	con/.style      = {treenode, font=\ttfamily, bottom color=green!25},
	nocon/.style    = {treenode, font=\ttfamily, bottom color=red!30},
	dummy/.style    = {circle,draw}
}

\startlocaldefs

\newcommand{\eat}[1]{}

\usepackage{enumitem}
\DeclareMathOperator*{\argmin}{\arg\!\min}

\newcommand{\hk}{\mathcal{H}_K}
\newcommand{\kmac}{\hat{\eta}_{n}}
\renewcommand{\hat}[1]{\widehat{#1}}
\renewcommand{\tilde}[1]{\widetilde{#1}}
\newcommand{\E}{\mathbbm{E}}
\renewcommand{\P}{\mathbbm{P}}
\newcommand{\X}{{\mathcal{X}}}
\newcommand{\Y}{{\mathcal{Y}}}

\newcommand{\emgn}{\mathcal{E}(\mathcal{G}_n)}

\newcommand{\R}{\mathbb{R}}

\newcommand\independent{\protect\mathpalette{\protect\independenT}{\perp}}
\def\independenT#1#2{\mathrel{\rlap{$#1#2$}\mkern2mu{#1#2}}}

\newcommand{\mcg}{\mathcal{G}}

\newcommand{\tmk}{{\mathcal{M}}_K}

\@addtoreset{proofpart}{theorem}

\makeatletter
\newcommand*{\rom}[1]{\expandafter\@slowromancap\romannumeral #1@}
\makeatother

\def\argmin{\mathop{\rm argmin}}

\newcommand{\rgn}{\mathcal{G}_n^{\textsc{rank}}}
\newcommand{\rmgn}{\mathcal{E}(\mathcal{G}_n^{\textsc{rank}})}

\newcommand{\trgn}{\tilde{\mathcal{G}}_n^{\textsc{rank}}}
\newcommand{\rkmac}{\hat{\eta}_n^{\textsc{rank}}}
\newcommand{\prkmac}{\eta_K^{\textsc{rank}}}

\newcommand{\hS}{\tilde{S}}
\newcommand{\td}{\tilde{d}}

\newcommand{\nrk}{N_n^{\textsc{rank}}}
\newcommand{\nrp}{N_n^{\textsc{pop}}}

\newcommand{\mme}{\mathbbm{E}}

\newcommand{\rpk}{\eta_K^{\textsc{rank}}(\mu)}

\definecolor{LightCyan}{rgb}{0.88,1,1}
\definecolor{Gray}{gray}{0.9}

\endlocaldefs

\usepackage{abstract}
\usepackage{xr}
\RequirePackage[sort,numbers]{natbib}
\begin{document}
	
	\renewcommand{\abstractname}{}    
	\renewcommand{\absnamepos}{empty}
	
	\begin{frontmatter}
		\title{Distribution-free  Measures of Association   based on Optimal Transport}
		
		\runtitle{Distribution-free Measures of Association}

		\begin{aug}
			\author[A]{\fnms{Nabarun}~\snm{Deb}\ead[label=e1]{nabarun.deb@chicagobooth.edu}},
			\author[B]{\fnms{Promit}~\snm{Ghosal}\thanksref{t2}\ead[label=e2]{promit@uchicago.edu}}
			\and
			\author[C]{\fnms{Bodhisattva}~\snm{Sen}\thanksref{t1}\ead[label=e3]{bodhi@stat.columbia.edu}}
			
			\thankstext{t1}{Supported by NSF grant DMS-2311062.}

                \thankstext{t2}{Supported by NSF grant DMS-2346685.}
			
			\runauthor{Deb, Ghosal, and Sen}
			
			\address[A]{University of Chicago\printead[presep={,\ }]{e1}}
			
			\address[B]{University of Chicago\printead[presep={,\ }]{e2}}

                \address[C]{Columbia University\printead[presep={,\ }]{e3}}
			
			\end{aug}
		\vspace{0.2in}
		\begin{abstract}
			In this paper we propose and study a class of nonparametric, yet interpretable measures of association between two random vectors $X$ and $Y$ taking values in $\R^{d_1}$ and $\R^{d_2}$ respectively ($d_1, d_2\ge 1$). These nonparametric measures --- defined using the theory of reproducing kernel Hilbert spaces coupled with optimal transport --- capture the strength of dependence between $X$ and $Y$ and have the property that they are 0 if and only if the variables are independent and 1 if and only if one variable is a measurable function of the other. Further, these population measures can be consistently estimated using the general framework of geometric graphs which include $k$-nearest neighbor graphs and minimum spanning trees. Additionally, these measures can also be readily used to construct an exact finite sample distribution-free test of mutual independence between $X$ and $Y$. 
		   In fact, as far as we are aware, these are the only procedures that possess all the above mentioned desirable properties. The correlation coefficient proposed in~\citet{dette2013copula},~\citet{chatterjee2019new} and~\citet{Azadkia-Chatterjee-2021}, at the population level, can be seen as a special case of this general class of measures.
		\end{abstract}

		
		
	\end{frontmatter}

\section{Introduction}
Suppose that $Z=(X,Y)\sim\mu$ where $\mu$ is an absolutely continuous distribution (with respect to the Lebesgue measure) supported on some subset of $\R^{d_1+d_2}$ ($d_1,d_2 \ge 1$) with marginal distributions $\mu_X$ and $\mu_Y$ on $\X \subset \R^{d_1}$ and $\Y \subset \R^{d_2}$ respectively. Assume that we have i.i.d.~data $\{Z_i \equiv ({X}_i,{Y}_i)\}_{i = 1}^{n}$ from $\mu$. Our goal is to construct a real-valued measure to quantify the degree of association or dependence between $X$ and $Y$, both in the population and sample settings. 

Pearson's correlation coefficient (see e.g.,~\cite{pearson1920notes}) is perhaps the simplest measure of association between $X$ and $Y$ when $d_1 = d_2 =1$.  However, the classical correlation has the following drawback: It is not effective for detecting associations that are nonlinear in nature, even in the complete absence of noise. Recently,~\citet{chatterjee2019new} (also see~\cite{dette2013copula}) proposed a {\it nonparametric} measure of association between $X$ and $Y$ when $d_1 = d_2 =1$ which is 0 if and only if $X$ and $Y$ are independent and 1 if and only if $Y$ is a measurable function of $X$. Moreover, any value between 0 and 1 of the coefficient conveys an idea of the strength of the relationship between $X$ and $Y$. In~\citet{Azadkia-Chatterjee-2021} the authors extend this measure to the case when $d_1 \ge 1$ and $d_2 =1$.~\citet{DebEtAl-2020} propose a class of simple, nonparametric, yet interpretable measures of association between $X$ and $Y$ when they take values in general metric spaces. The above papers also propose and study consistent estimators of these population measures of dependence  using ideas based on {\it geometric graphs} (e.g., $k$-nearest neighbor graph, minimum spanning tree, etc.).

To motivate our contribution in this paper consider the Spearman's correlation between $X$ and $Y$ (see~\cite{spearman1904proof}) which is another simple and useful measure of dependence when $d_1 = d_2 = 1$. An important property of the sample Spearman's correlation coefficient (not possessed by Pearson's correlation coefficient) is that it is {\it distribution-free}, i.e., its distribution does not depend on $\mu$ when $X$ and $Y$ are independent  (also see Kendall's $\tau$~\cite{kendall1938new,Kendall1990}). This, in particular, implies that the sample Spearman's correlation coefficient can not only be used as a measure of dependence but can also be readily used to construct an {\it exact distribution-free procedure} to test the hypothesis of mutual independence between $X$ and $Y$ (when $d_1 = d_2 = 1$).

This leads to the following question: Can we construct a \textit{nonparametric measure of association} between $Y$ and $X$ such that its sample analogue is {\it distribution-free} when $d_1, d_2 \ge 1$? 
    

In this paper, we answer this question in the affirmative and propose and study a class of {\it distribution-free} empirical measures $T_n\equiv T_n(Z_1,\ldots ,Z_n)$ and their population counterparts when $d_1, d_2 \ge 1$, that yield a family of {\it nonparametric measures of association} between ${Y}$ and ${X}$. Our work extends the framework of~\citet{DebEtAl-2020}. In particular, these measures are defined using the theory of reproducing kernel Hilbert spaces (RKHS); see Section~\ref{sec:rkhsprelim} for a brief review. 

A plethora of nonparametric procedures have been proposed that can detect nonlinear dependencies between the variables $X$ and $Y$ over the last 60 years; see e.g.,~\cite{Renyi59, blum1961distribution, Rosenblatt75, friedman1983, Gabor2007, gieser1997, taskinen2005, gretton2008kernel, Oja-2010,Reshef11,Lyons13, heller2013, SS-2014-Bio, JH16, munmun2016, berrett2017nonparametric} and the references therein. While these coefficients are useful in practice, they have one common problem: They are all designed primarily for testing independence, and not for measuring the strength of the relationship between the variables. Moreover, as far as we are aware, when $d_1, d_2 > 1$ none of above sample measures are distribution-free when $X$ and $Y$ are independent.

To construct our distribution-free measures of association we borrow ideas from three seemingly disparate fields: (i) The theory of RKHS which is needed to define the class of {\it nonparametric} measures of association; (ii) the notion of geometric graphs (see Section~\ref{sec:gengraphcond} for a review) which are used to develop {\it consistent} estimators of these nonparametric measures; and (iii) the theory of {\it optimal transport} (OT) which yields {\it distribution-free} empirical measures.


 


Our starting point is the class of kernel based measures of association defined in~\citet{DebEtAl-2020} (see Section~\ref{sec:Kmac} for a brief review). Although the proposed empirical measure there (denoted there by $\kmac$) is nonparametric, it is not distribution-free. A key observation in this regard is that the distribution-free measures discussed above when $d_1=d_2 = 1$ (e.g., Spearman's correlation coefficient) are based on the (univariate) ``ranks'' of $Y_i$'s. Hence, one could imagine that a `proper' distribution-free notion of multivariate ranks can be employed to construct distribution-free notions of dependence when $d_1, d_2 >1$. Indeed this is our approach: We use the recently developed idea of {\it multivariate ranks} based on the theory of OT (see Sections~\ref{sec:overviewOT} and~\ref{sec:defn} for a brief review on this topic; also see~\cite{Cher17,Deb19,Shi19}) to develop measures of association that are finite-sample distribution-free when $\mu=\mu_X\otimes\mu_Y$.

Having defined the multivariate ranks (via OT) we construct our family of distribution-free measures of association based on a very simple and classical analogy between Pearson's correlation and Spearman's correlation. Note that when $d_1=d_2=1$, Spearman's correlation is equivalent to the classical Pearson's correlation coefficient computed between the one-dimensional ranks of the $X_i$'s and the $Y_i$'s, instead of using the observations themselves. We mimic the same approach here, i.e., instead of computing $\kmac$ (as in~\cite{DebEtAl-2020}) using the $X_i$'s and $Y_i$'s themselves, we instead use their empirical multivariate ranks. 
	
We propose the ``rank" version of $\kmac$, namely $\rkmac$, in~\eqref{eq:rankkmac}. In~\cref{thm:Consistency}, we show that $\rkmac$ consistently estimates a population measure of dependency $\prkmac$ between $X$ and $Y$, and is distribution-free when $\mu=\mu_X\otimes\mu_Y$. In~\cref{thm:rankassoc} we show that this population measure $\prkmac$ satisfies many desirable properties that justifies it as a population measure of association between $X$ and $Y$. In particular, we show that: 
\begin{itemize}
		\item[(P1)] $\prkmac  \in [0,1]$;
		\item[(P2)] $\prkmac=0$ if and only if $\mu=\mu_X\otimes \mu_Y$ (i.e., $X$ and $Y$ are independent);
		\item[(P3)] $\prkmac=1$ if and only if $Y=g(X)$, $\mu$ almost everywhere (a.e.), for some measurable function $g:\R^{d_1} \to \R^{d_2}$.
	\end{itemize}

As $\rkmac$ is based on multivariate ranks, a test for independence of $X$ and $Y$ based on $\rkmac$ will generally be more powerful against heavy-tailed alternatives and more robust to outliers and contamination (see~\cite{Huber-2nd-E, Oja-2010, Deb19} for related discussions). Further, the corresponding test, being distribution-free, also avoids asymptotic approximations or permutation ideas for determining rejection thresholds.

In~\cref{prop:Chacon}, we prove that the limit of $\rkmac$, i.e., $\prkmac$  exactly coincides with the limit of the coefficient in~\citet{Azadkia-Chatterjee-2021} (denoted by $T_n(Y,Z)$ in their paper) when $d_2=1$ for an appropriate choice of a kernel. Note that, unlike $\rkmac$, the empirical  measure in~\cite{Azadkia-Chatterjee-2021} does not have the finite sample distribution-free property.
	
Finally,~\cref{theo:ranknullclt} proves a central limit theorem for $\rkmac$ which is uniform over a large class of (geometric) graphs. We would like to point out that unlike $X_i$'s and $Y_i$'s, their multivariate ranks are no longer independent among themselves which makes the CLT challenging to prove.

The paper is organized as follows: In Section~\ref{sec:Prelims} we briefly review the RKHS framework, the notion of geometric graphs, and multivariate ranks via OT. Our distribution-free nonparametric measure of association is introduced in Section~\ref{sec:R-Q-Intro} along with the main results of this paper. We conclude with some remarks and open questions in Section~\ref {sec:Disc}. Section~\ref{sec:pfmain} gives the proofs of all the main results in this paper. 
	
\subsection{Related works} 
In~\citet{dette2013copula}, the authors use the term ``measure of regression dependence" for the three properties mentioned above and show that it is possible to define such a measure satisfying (P1)-(P3) when $\mathcal{X}=\Y=\R$. The same population measure was rediscovered in~\citet{chatterjee2019new} where the author also proposed a tuning parameter-free estimator of the same measure that can be computed in near linear time. Since then, the estimator in~\cite{chatterjee2019new} has attracted a lot of attention (see e.g.,~\cite{Shi-Drton-Han-2022,cao2020correlations}). Further, in~\citet{Azadkia-Chatterjee-2021}, the authors propose a similar measure when $\mathcal{X}=\R$ and $\Y=\R^{d_2}$, $d_2\geq 1$. However, all these measures crucially use the canonical ordering of $\R$ and hence do not extend to the multivariate setting (where $\mathcal{X}=\R^{d_1}$ and $\Y=\R^{d_2}$ with $d_1,d_2 \ge 2$). Some multivariate measures of association satisfying (P1)-(P3) have been proposed in~\cite{siburg2010measure,tasena2016measures,boonmee2016measure}, following similar copula-based ideas as in~\cite{dette2013copula}; however, to the best of our knowledge, neither of these papers provide a consistent empirical estimate of their proposed measures of association, nor do the associated tests for independence have the distribution-free property. 
 
 In a number of other papers, e.g.,~\cite{Shi-Drton-Han-2022, nies2021transport, Auddy-Et-Al-2024, dette2023simple,zhang2023asymptotic,lin2023boosting}, the authors analyze different properties of the estimator in \citet{chatterjee2019new}. In \cite{ansari2022}, the authors provide an extension of the measure in \cite{Azadkia-Chatterjee-2021} to a vector of endogenous variables whereas \cite{fuchs2023quantifying} connects copula based measures of association to dimension reduction principles. In another line of work \cite{strothmann2024rearranged}, the authors provide a recipe to modify a given measure of dependence to yield a measure of association satisfying (P1)-(P3) above. Finally, we refer the reader to~\cite{chatterjee2022survey,chattamvelli2024measures} for a survey of the recent papers in this topic.


\section{Preliminaries}\label{sec:Prelims}
In this section we briefly introduce three topics: (i) The (nonparametric) class of kernel measures of association as defined in~\citet{DebEtAl-2020} using the theory of RKHS; (ii) the notion of geometric graphs which will be used to develop consistent estimators of these nonparametric population measures; and (iii) the theory of OT  which will be used crucially to obtain distribution-free empirical measures.

Let $\X$ and $\Y$ be topological spaces. Let $(X,Y) \in \X \times \Y$ with joint distribution $\mu$ and marginals $X \sim \mu_X$ and $Y \sim \mu_Y$. Although we will eventually take $\X = \R^{d_1}$ and $\Y = \R^{d_2}$, for much of what we review in this section we do not actually need to make this assumption. Let $\mathcal{P}(\Y)$ and $\mathcal{P}(\X \times \Y)$ be the set of all Borel probability measures on $\Y$  and $\X \times \Y$ respectively. Suppose that $\mu$ admits a regular conditional distribution $\mu_{Y|x}$ --- the conditional distribution of $Y$ given $X = x$; existence of regular conditional distributions can be guaranteed under mild conditions (see~\cite{Faden1985} for a survey). 

\subsection{RKHS: Some background}\label{sec:rkhsprelim}
In this subsection we formally define some concepts from the theory of RKHS that will be used repeatedly in this paper. We start with the basic definition of a RKHS.
	
\begin{defn}[Reproducing kernel Hilbert space (RKHS)]\label{def:rkhs}
Let $\mathcal{H}$ be a Hilbert space of real-valued functions defined on $\Y$ with inner product $\langle \cdot,\cdot\rangle_{\mathcal{H}}: \mathcal{H} \times \mathcal{H} \to \R$. A function $K:\Y\times\Y\to\R$ is called a \emph{reproducing kernel} if the following two conditions hold:
		\begin{enumerate}
			\item For all $y\in\Y$, $K(\cdot,y)\in\mathcal{H}$.
			\item For all $y\in\Y$ and $f\in\mathcal{H}$, $\langle f,K(\cdot,y)\rangle_{\mathcal{H}}=f(y)$.
		\end{enumerate} 
		If $\mathcal{H}$ admits a reproducing kernel, then it is termed as a RKHS.
	\end{defn}
By the Moore-Aronszajn Theorem (see e.g.,~\cite[Theorem 3]{Berlinet2004}) a symmetric, nonnegative definite kernel function $K(\cdot,\cdot)$ on $\Y\times \Y$ can be identified uniquely with a unique RKHS of real-valued functions on $\Y$ for which $K(\cdot,\cdot)$ is the reproducing kernel. Let us denote this RKHS by $\mathcal{H}_K$. Let us also assume that $\hk$ is separable (this can be ensured under mild conditions\footnote{For example, if $\Y$ is a separable space and $K(\cdot,\cdot)$ is continuous.}, see e.g.,~\cite[Lemma 4.33]{SVM}).
	
The map $y\mapsto K(\cdot,y)$ from $\Y$ to $\hk$ is often called the \emph{feature map}. Further, the reproducing property as stated in~\cref{def:rkhs} implies that 
	\begin{equation}\label{eq:hilnorm}
	\langle K(\cdot,y),K(\cdot,\tilde{y})\rangle_{\hk}=K(y,\tilde{y}), \qquad \mbox{for all}\;y,\tilde{y} \in \Y.
	\end{equation}	
In the following we define three concepts that will be crucial in defining $\eta_K$ --- the kernel measure of dependence as introduced in~\citet{DebEtAl-2020}. Suppose that $Y \sim \mu_Y$ has a probability distribution on $\Y$ such that $\E[\sqrt{K(Y,Y)}] < \infty$. 
	
	\begin{defn}[Mean embedding]\label{def:meanembed}
		Define the following class of probability measures on $\Y$: $$\tmk^{\theta}(\Y)\coloneqq \left\{\nu\in\mathcal{P}(\Y):\int_\Y K^{\theta}(y,y)\,d\nu(y)<\infty \right\}, \qquad \mbox{for } \, \theta>0.$$ Let $\mu_Y\in\tmk^{1/2}(\Y)$. Then the (kernel) \emph{mean embedding} of $\mu_Y$ into $\hk$ is given by $m_K(\mu_Y)\in \hk$ such that 
		\begin{equation}\label{eq:kmeanembed}
		\langle f,m_K(\mu_Y)\rangle_{\hk} = \int_\Y f(y)\,d\mu_Y(y), \qquad \mbox{for all } \, f\in\hk.
		\end{equation}
	\end{defn}
	In fact, one can write $m_K(\mu_Y)=\int K(\cdot,y)\,d\mu_Y(y)=\E_{\mu_Y} [K(\cdot,Y)]$. It is well-defined as a consequence of the Riesz representation theorem,~see~\cite{smola2007hilbert,akhiezer1993} (equivalently also by Bochner's theorem, see~\cite{Diestel1974,Dinculeanu2011}). The map $\mu_Y\mapsto m_K(\mu_Y)$ with domain $\tmk^{1/2}(\Y)$ can be viewed as a natural extension of the map $y\mapsto K(\cdot,y)$ with domain $\Y$.
	
	In a similar vein we can also define the (kernel) \emph{mean embeddings} of conditional distributions. For $X=x$, the (kernel) \emph{conditional mean embedding} of $\mu_{Y|x}$ is defined as an element of $\hk$ in the same way as~\eqref{eq:kmeanembed} (also see~\cite[Section 4.1.1]{muandet2017kernel}). In other words, we can also write $m_K(\mu_{Y|x})=\int K(\cdot,y)d\mu_{Y|x}(y)=\E_{\mu_{Y|x}}[K(\cdot,Y)]$.

	\begin{defn}[Maximum mean discrepancy]\label{defn:MMD}
		The difference between two probability distributions $Q_1$ and $Q_2$ in $\tmk^{1}(\Y)$ can then be conveniently measured by $${\rm MMD}_K(Q_1,Q_2) := \|m_K(Q_1) - m_K(Q_2)\|_{\hk}$$ (here $m_K(Q_i)$ is the mean element of $Q_i$, for $i=1,2$) which is called the {\it maximum mean discrepancy} (MMD) between $Q_1$ and $Q_2$ (see~\cite[Definition 10]{gretton2008}). The following alternative representation of the squared MMD is also known (see e.g.,~\cite[Lemma 6]{Gretton12}, or simply use~\eqref{eq:hilnorm}):
		\begin{equation}\label{eq:MMD}
		{\rm MMD}_K^2(Q_1,Q_2) = \E[K(S,S')] + \E[K(W,W')] - 2 \E[K(S,W)],
		\end{equation}
		where $S, S',W,W'$ are independent, $S, S' \stackrel{}{\sim} Q_1$ and $W, W' \stackrel{}{\sim} Q_2$.
	\end{defn}
	
	\begin{defn}[Characteristic kernel]\label{defn:Char}
		The kernel $K(\cdot,\cdot)$ is said to be {\it characteristic} if and only if the map $Q\mapsto m_K(Q)$ is one-to-one on the domain $\tmk^1(\Y)$, i.e., $$m_K(Q_1) = m_K(Q_2) \quad  \implies \quad Q_1=Q_2, \quad \mbox{for all }\;Q_1,Q_2\in \tmk^1(\Y).$$
	\end{defn}
	Note that the last condition is equivalent to $\langle m_K(Q_1),f\rangle_{\hk} = \langle m_K(Q_2),f\rangle_{\hk}$ for all $f \in \hk$, i.e., $\E_{S \sim Q_1}[f(S)] = \E_{W \sim Q_2}[f(W)]$, for all $f \in \hk$. A characteristic kernel implicitly implies that the associated RKHS is rich enough. 

\subsection{A kernel measure of association}\label{sec:Kmac}
In this subsection we review the class of nonparametric kernel measures of association introduced in~\citet{DebEtAl-2020}. 
Recall that $\mu_{Y|x}$ is the regular conditional distribution of $Y$ given $X=x$ which we assume exists for all $x \in \X$. Assume that a {\it kernel} $K(\cdot,\cdot)$  --- a symmetric, nonnegative definite function on $\Y\times \Y$ --- exists on $\Y\times\Y$ and let $\hk$ denote the induced separable RKHS (see~\cref{def:rkhs}). Let $\langle \cdot,\cdot\rangle_{\hk}: \hk \times \hk\to \R$ and $\|\cdot \|_{\hk}$ denote the inner product and induced norm on $\hk$.

Generate $(X',Y',\tilde{Y'})$ as follows: $X'\sim\mu_X$, $Y'|X\sim \mu_{Y|X'}$, $\tilde{Y'}|X'\sim \mu_{Y|X'}$ and $Y'\independent\tilde{Y'}|X'$. Also let $Y_1,Y_2$ be i.i.d. from $\mu_Y$. Note that $Y'$ and $\tilde{Y'}$ are dependent (via $X'$), unlike $Y_1$ and $Y_2$. The measure of association of $Y$ on $X$ as introduced in~\citet{DebEtAl-2020} can now be presented as:
	\begin{equation}\label{eq:kac}
	\eta_{K}(\mu) := 1-\frac{\mme\lVert K(\cdot,Y')-K(\cdot,\tilde{Y'})\rVert_{\mathcal{H}_K}^2}{\mme\lVert K(\cdot,Y_1)-K(\cdot,Y_2)\rVert_{\mathcal{H}_K}^2}.
	\end{equation}
In order to ensure that $\eta_{K}(\cdot)$ is well-defined, we need certain moment assumptions. By the reproducing property of $K(\cdot,\cdot)$,
	\begin{equation}\label{eq:K-1-2}	
	\lVert K(\cdot,Y_1)-K(\cdot,Y_2)\rVert_{\mathcal{H}_{K}}^2=K(Y_1,Y_1)+K(Y_2,Y_2)-2K(Y_1,Y_2). 
	\end{equation}
Suppose that $\mu_Y\in\tmk^1(\Y)$.
Then the first two terms in~\eqref{eq:K-1-2} have finite moments. The third term is also finite by an application of the Cauchy-Schwartz inequality combined with the observation that $\tmk^1(\Y)\subset \tmk^{1/2}(\Y)$. Thus, we assume that $\mu_Y\in \tmk^1(\Y)$ in the sequel. In the following we give an alternate expression of $\eta_K$ which will be crucial and can be easily derived from~\eqref{eq:kac} and~\eqref{eq:K-1-2} (also see~\cite{DebEtAl-2020}):
\begin{equation}\label{eq:Alt-eta-K}
\eta_K(\mu) = \frac{\E[ K(Y_1,Y_2)]-\E [K(Y',\tilde{Y'})]}{\E [K(Y_1,Y_1)]-\E[ K(Y_1,Y_2)]}.
\end{equation}

The following result shows that $\eta_K$ is indeed a valid measure of association with many desirable properties. 
	\begin{theorem}[Theorem 2.1 of~\cite{DebEtAl-2020}]\label{theo:kacmas}
		Suppose $\mu_Y\in \tmk^1(\Y)$, $\Y$ is Hausdorff\,\footnote{A topological space where for any two distinct points there exist neighborhoods of each which are disjoint from each other.} and $K(\cdot,\cdot)$ is a characteristic kernel. Then $\eta_K(\mu)$, as defined in~\eqref{eq:kac}, satisfies the following properties: 	\begin{itemize}
		\item[(P1)] $\eta_K(\mu) \in [0,1]$.
		\item[(P2)] $\eta_K(\mu)=0$ if and only if $\mu=\mu_X\otimes \mu_Y$ (i.e., $X$ and $Y$ are independent).
		\item[(P3)] $\eta_K(\mu)=1$ if and only if $Y=g(X)$, $\mu$ almost everywhere (a.e.), for some measurable function $g:\mathcal{X}\to\Y$.
	\end{itemize}

	\end{theorem}

\subsection{Estimation via geometric graphs}\label{sec:gengraphcond}
In this subsection we review estimation of $\eta_K$ (as in~\cite{DebEtAl-2020}) using ideas from geometric graphs. Note that the denominator in~\eqref{eq:Alt-eta-K} can be easily estimated using empirical averages (from standard U-statistics theory; see~\cite[Chapter 12]{van1998}), for instance, with the following estimator: 
$$\frac{1}{n} \sum_{i=1}^n K(Y_i,Y_i) - \frac{1}{n(n-1)}\sum_{i\neq j} K(Y_i,Y_j).$$
The numerator in~\eqref{eq:Alt-eta-K} is trickier to estimate. The main difficulty arises because of the term $\E\big[K(Y',\tilde{Y'})\big] = \E\big[\E(K(Y',\tilde{Y'}) \mid X')\big]$ as we do not observe two $Y$ values --- $Y'$ and $\tilde{Y'}$ --- independently from the conditional distribution of $X'$ (cf.~the notation introduced just before~\eqref{eq:kac}). To overcome this difficulty, to estimate the above term we consider a `neighbor' of $X'$, say $\tilde{X'}$, and its corresponding $Y$ value, which can be used as a surrogate to $\tilde{Y'}$. We formalize this notion of a `neighbor' using a geometric graph on $\X$. 

$\mcg$ is a \emph{geometric graph functional} on $\mathcal{X}$ if, given any finite subset $S \subset \mathcal{X}$, $\mcg(S)$ defines a graph with vertex set $S$ and corresponding edge set, say $\mathcal{E}(\mcg(S))$.  
Note that the graph $\mathcal{G}(S)$ can be \emph{directed/undirected}. In this paper, we will restrict ourselves to simple graphs (i.e., those without multiple edges and self loops) with no isolated vertices. Accordingly, we will often drop the qualifier geometric and simple. 
	
Next we define $\mathcal{G}_n\coloneqq \mathcal{G}(X_1,\ldots ,X_n)$ where $\mathcal{G}$ is some graph functional on $\mathcal{X}$. We would like to define graph functionals for which $(i,j)\in\emgn$ implies $X_i$ and $X_j$ are ``close".~\citet{DebEtAl-2020} considered the following sample analogue of $\E\big[K(Y',\tilde{Y'})\big]$:
\begin{equation}\label{eq:analogtarget}
	\frac{1}{n}\sum_{i=1}^n\frac{1}{d_i}\sum_{j:(i,j)\in\emgn} K(Y_i,Y_j),
\end{equation}
where $\emgn$ denotes the set of (directed/undirected) edges of $\mathcal{G}_n$, i.e., $(i,j)\in\emgn$ if and only if there is an edge from $i\to j$ or $j\to i$ in $\mathcal{G}_n$, and $d_i$ denotes the degree of $X_i$ in $\mathcal{G}_n$. To be specific $d_i:=\sum_{j:(i,j)\in\emgn} 1$. 
	
Using~\eqref{eq:analogtarget},~\citet{DebEtAl-2020} proposed the following estimator of $\eta_K$:
	\begin{align}\label{eq:statest}
	\kmac:=\frac{\frac{1}{n}\sum_{i=1}^n d_i^{-1}\sum_{j:(i,j)\in\emgn} K(Y_i,Y_j)- \frac{1}{n(n-1)}\sum_{i\neq j} K(Y_i,Y_j)}{\frac{1}{n} \sum_{i=1}^n K(Y_i,Y_i) - \frac{1}{n(n-1)}\sum_{i\neq j} K(Y_i,Y_j)}.
	\end{align}
	
	
Under suitable conditions on the graph functional,~\cite{DebEtAl-2020} showed that indeed $\kmac$ consistently estimate $\eta_K$ as the sample size grows to infinity.

\subsection{A brief overview of OT}\label{sec:overviewOT}
In this subsection we introduce the basics of OT which will be necessary to construct distribution-free estimates of our measures of association. Let $\mathcal{P}_{ac}(\R^d)$ denote the space of absolutely continuous probability measures on $\R^d$. For a function $F:\R^d\to\R^d$, we will use $F\#\mu$ to denote the push forward measure of $\mu$ under $F$, i.e., the distribution of $F(Z)$ when $Z\sim\mu$.
	
Below we present perhaps the simplest version of the OT problem (courtesy the works of Gaspard Monge in 1781, see~\cite{monge1781memoire}):
	\begin{align}\label{eq:Mongeproblem}
	\inf_{F} \int  \lVert z-F(z)\rVert^2 \,d\mu(z)\qquad \mbox{subject to}\quad F\#\mu=\nu.
	\end{align}
A minimizer of~\eqref{eq:Mongeproblem}, if it exists, is referred to as an {\it OT map}. An important result in this field, known as the {\it Brenier-McCann theorem}, takes a ``geometric" approach to the problem of OT (as opposed to the analytical approach presented in~\eqref{eq:Mongeproblem}) and will be very useful to us in the sequel; see e.g.,~\cite[Theorem 2.12 and Corollary 2.30]{Villani2003}.

\begin{prop}[Brenier-McCann  theorem~\cite{Mccann1995}]\label{prop:Mccan}
		Suppose that $\mu, \nu\in\mathcal{P}_{ac}(\R^d)$. Then there exists functions $R(\cdot)$ and $Q(\cdot)$ (usually referred to as ``OT maps"), both of which are gradients of (extended) real-valued $d$-variate convex functions, such that: (i) $R\# \mu=\nu$, $Q\# \nu=\mu$; (ii) $R$ and $Q$ are unique ($\mu$ and $\nu$-a.e.~respectively); (iii) $R\circ Q (u)=u$ ($\nu$-a.e.~$u$) and $Q\circ R(z)=z$ ($\mu$-a.e.~$z$). Moreover, if $\mu$ and $\nu$ have finite second moments, $R(\cdot)$ is also the solution to Monge's problem in~\eqref{eq:Mongeproblem}.
\end{prop}
In~\cref{prop:Mccan}, by ``gradient of a convex function" we essentially mean a function from $\mathbb{R}^d$ to $\mathbb{R}^d$ which is $\mu$ (or $\nu$) a.e.~equal to the gradient of some convex function. It is instructive to note that, when $d=1$, the standard $1$-dimensional distribution function $F$ associated with a distribution $\mu$ is nondecreasing and hence the gradient of a convex function. Therefore when $d=1$, $F$ is the OT map from $\mu$ to $\mathcal{U}[0,1]$ (i.e., the Uniform([0,1]) distribution) by~\cref{prop:Mccan}.
	
	\subsection{Multivariate ranks defined via OT}\label{sec:defn}
	\begin{defn}[Population multivariate ranks and quantiles]\label{def:popquanrank}
		Set $\nu:=\mathcal{U}[0,1]^d$ --- the uniform distribution on $[0,1]^d$. Given $\mu\in\mathcal{P}_{ac}(\mathbb{R}^d)$, the corresponding population rank and quantile maps are defined as the OT maps $R(\cdot)$ and $Q(\cdot)$ respectively as in~\cref{prop:Mccan}. These maps are unique a.e.~with respect to $\mu$ and $\nu$ respectively. 
	\end{defn}
	
	In standard statistical applications, the population rank map is not available to the practitioner. In fact, the only accessible information about $\mu$ comes in the form of empirical observations $Z_1,Z_2,\ldots ,Z_n\overset{i.i.d.}{\sim}\mu\in\mathcal{P}_{ac}(\mathbb{R}^d)$. In order to estimate the population rank map from these observations, let us denote
	\begin{equation}\label{eq:H_n^d}
	\mathbf{\mathcal{H}}_n^d:=\{h_1^d,\ldots ,h_n^d\}
	\end{equation}
	to be a set of $n$ vectors in $[0,1]^d$. We would like the points in $\mathbf{\mathcal{H}}_n^d$ to be ``uniform-like", i.e., their empirical distribution, $n^{-1}\sum_{i=1}^n \delta_{h_i^d}$ should converge weakly to $\nu$. In practice, for $d=1$, we may take $\mathbf{\mathcal{H}}_n^d$ to be the usual $\{i/n\}_{1\leq i\leq n}$ sequence and for $d\geq 2$, we may take it as a quasi-Monte Carlo sequence (such as the $d$-dimensional Halton or Sobol sequence) of size $n$ (see~\cite{Hofer2009,Hofer2010} for details) or a random draw of $n$ i.i.d. random variables from $\nu$. The empirical distribution on $\mathbf{\mathcal{H}}_n^d$ will serve as a discrete approximation of $\nu$. We are now in a position to define the empirical multivariate rank function which will proceed via a discrete analogue of problem~\eqref{eq:Mongeproblem}.
	
	\begin{defn}[Empirical rank map]\label{def:empquanrank}
		Let $S_n$ denote the set of all $n!$ permutations of $\{1,2,\ldots ,n\}$. Consider the following optimization problem:  
		\begin{align}\label{eq:empopt}
		\mathbf{\hat{\sigma}}_n\coloneqq \argmin_{\mathbf{\sigma}\in S_n} \sum_{i=1}^n \lVert X_{i}-h_{{\sigma(i)}}^d\rVert^2.
		\end{align}  
		Note that $\mathbf{\hat{\sigma}}_n$ is a.s.~uniquely defined (for each $n$) as $\mu\in\mathcal{P}_{ac}(\R^d)$. The empirical ranks are then defined as:
		\begin{equation}\label{eq:R_n_h}
		\hat{R}_n(X_i)=h_{\hat{\sigma}_n(i)}^d, \qquad \mathrm{for \;\;} i = 1,\ldots, n.
		\end{equation}
	\end{defn}
	The optimization problem in~\eqref{eq:empopt} is combinatorial in nature, but it can be solved exactly in polynomial time (with worst case complexity $\mathcal{O}(n^3)$) using the Hungarian algorithm (see~\cite{munkres1957,bertsekas1988} for details). For a comprehensive list of faster (approximate) algorithms, see~\cite[Section 5]{Shi19}. Moreover, when $d=1$, if we choose $\mathbf{\mathcal{H}}_n^1= \{i/n\}_{i=1}^{n}$, then the empirical ranks, i.e., $\hat{R}_n(X_i)$'s are exactly equal to the usual one-dimensional ranks.
	
	A crucial reason behind defining empirical ranks as in~\cref{def:empquanrank} is that due to the exchangeability of the $X_i$'s, the vector of ranks, i.e., $(\hat{R}_n(X_1),\ldots ,\hat{R}_n(X_n))$ is uniformly distributed over the following set:
	$$\{(h^d_{\sigma(1)},\ldots ,h^d_{\sigma(n)}):\sigma\in S_n\}.$$
	This lends a distribution-free property to the empirical ranks (see~\cite[Proposition 2.2]{Deb19} for a formal statement). In other words, the distribution of $(\hat{R}_n(X_1),\ldots ,\hat{R}_n(X_n))$ is free of $\mu\in\mathcal{P}_{ac}(\R^d)$. Moreover, the empirical ranks are maximal ancillary (see e.g.,~\cite{del2018center}).
	

\section{A multivariate rank based measure of association}\label{sec:R-Q-Intro}
In this section we develop a new class of measures of association and their sample analogues which have the distribution-free property. We will restrict our attention to the case when $\X=\R^{d_1}$ and $\Y=\R^{d_2}$, for $d_1,d_2 \ge 1$.

\subsection{A distribution-free measure of association and its properties}
	
We will construct an empirical measure of association with the distribution-free property based on a simple and classical analogy between Pearson's and Spearman's correlation. Note that when $d_1=d_2=1$,  Spearman's correlation is equivalent to the classical Pearson's correlation coefficient computed between the one-dimensional ranks of the $X_i$'s and the $Y_i$'s, instead of the actual observations. It is this usage of one-dimensional ranks that lends Spearman's correlation the distribution-free property when $\mu=\mu_X\otimes\mu_Y$. We will mimic the same approach in this section, i.e., instead of computing $\kmac$ (as defined in~\eqref{eq:statest}) using the $X_i$'s and $Y_i$'s themselves, we will instead use their empirical multivariate ranks. 
	
Let us briefly recall the setting. We have $(X_1,Y_1),\ldots ,(X_n,Y_n)\overset{i.i.d.}{\sim}\mu\in \mathcal{P}_{ac}(\R^{d_1+d_2})$ with marginals $\mu_X$ and $\mu_Y$, such that $\mu_X\in\mathcal{P}_{ac}(\R^{d_1})$ and $\mu_Y\in \mathcal{P}_{ac}(\R^{d_2})$. Let $\mathbf{\mathcal{H}}_n^{d_1}$ and $\mathbf{\mathcal{H}}_n^{d_2}$ be two sets of $n$ ``uniform-like" points in dimension $d_1$ and $d_2$ respectively. We can then use~\cref{def:empquanrank} to define the empirical multivariate rank vectors $(\hat{R}_n^X(X_1),\ldots ,\hat{R}_n^X(X_n))$ and $(\hat{R}_n^Y(Y_1),\ldots ,\hat{R}_n^Y(Y_n))$ based on $\mathbf{\mathcal{H}}_n^{d_1}$ and $\mathbf{\mathcal{H}}_n^{d_2}$ (see~\eqref{eq:H_n^d}) respectively. 

Next, given a graph functional $\mathcal{G}$ (see  \cref{sec:gengraphcond} for its definition), let $$\rgn:=\mathcal{G}(\hat{R}_n^X(X_1),\ldots ,\hat{R}_n^X(X_n))$$ and $\rmgn$ be the set of edges of $\rgn$. Also, with a slight notational abuse, we will still use $(d_1,\ldots ,d_n)$ to be the degree sequence of the vertices in $\rgn$. The rank version of $\kmac$ is then defined as follows:
\begin{equation}\label{eq:rankkmac}
\rkmac:=\frac{n^{-1}\sum_{i}d_i^{-1}\sum_{j:(i,j)\in\rmgn} K(\hat{R}_n^Y(Y_i),\hat{R}_n^Y(Y_j))-F_n}{n^{-1}\sum_{i=1}^n K(\hat{R}_n^Y(Y_i),\hat{R}_n^Y(Y_i))-F_n},
	\end{equation}
	where 
\begin{equation}\label{eq:F_n}
F_n :=(n(n-1))^{-1}\sum_{i\neq j} K(\hat{R}_n^Y(Y_i),\hat{R}_n^Y(Y_j)).
\end{equation}
 Here $(i,j)\in\rmgn$ implies that $\hat{R}_n^Y(X_i)$ and $\hat{R}_n^Y(X_j)$ are connected in $\rgn$. 
	Note that $\rkmac$ is the same as $\kmac$ with the observations replaced with their empirical ranks. 
In the following theorems (proved in ~\cref{sec:pfmain}), we show that $\rkmac$ is a measure of association which additionally has a pivotal distribution when $\mu=\mu_X\otimes\mu_Y$.

\begin{theorem}\label{thm:Consistency}
		\emph{(a)} When $\mu=\mu_X\otimes\mu_Y$, $\rkmac$ has a pivotal distribution, i.e., its distribution does not depend on $\mu_X$ and $\mu_Y$.
		
		\emph{(b)} Assume that $K(\cdot,\cdot):\R^{d_2}\times\R^{d_2}\to\R$ is continuous and let $R^X(\cdot)$ and $R^Y(\cdot)$ denote the population rank maps (see~\cref{def:popquanrank}) for $\mu_X$ and $\mu_Y$ respectively. Suppose that the following assumptions hold:
        \begin{itemize}
        \item[(S1)] For $x_1,x_2 \in [0,1]^{d_1}$, define $$r(x_1,x_2):=\mathbb{E}\Big[K(R^Y(Y_1), R^Y(Y_2)) \mid R^X(X_1)=x_1, R^X(X_2)=x_2\Big]$$ and assume that it is uniformly $\beta$-H{\"o}lder continuous in $x_1,x_2\in [0,1]^{d_1}$ for some $\beta\in (0,1]$, i.e., given any $x_1,x_2,\tilde{x}_1,\tilde{x}_2\in [0,1]^{d_1}$, there exists a constant $C$ (free of $x_1,x_2,\tilde{x}_1,\tilde{x}_2$) such that:
		\begin{align}\label{eq:artifact}
		\big|r(x_1,x_2)-r(\tilde{x}_1,\tilde{x}_2)\big|\leq C\left(\lVert x_1-\tilde{x}_1\rVert^{\beta}+\lVert x_2-\tilde{x}_2\rVert^{\beta}\right).
		\end{align}
		\item[(S2)] Let $\rgn$ satisfy the following conditions: 
        \begin{equation}\label{eq:graph1}
        \limsup\limits_{n\to\infty}\frac{\max_{i=1}^n d_i}{\min_{i=1}^n d_i}<\infty,
        \end{equation}
        \begin{equation}\label{eq:graph2}
        \frac{1}{n\min_{i=1}^n d_i}\sum_{e\in \rmgn}|e|^{\beta} \to 0,   
        \end{equation}
        where $|e|$ denotes the edge length of $e\in\rmgn$.
        \item[(S3)] The empirical distributions on $\mathbf{\mathcal{H}}_n^{d_1}$ and $\mathbf{\mathcal{H}}_n^{d_2}$ converge weakly to $\mathcal{U}[0,1]^{d_1}$ and $\mathcal{U}[0,1]^{d_2}$ respectively. 
        \end{itemize}
        
        Under these assumptions, 
		$$\rkmac\overset{\mathbb{P}}{\longrightarrow} 1-\frac{\E\lVert K(\cdot,R^Y(Y'))-K(\cdot,R^Y(\tilde{Y'}))\rVert_{\hk}^2}{\E\lVert K(\cdot,R^Y(Y_1))-K(\cdot,R^Y(Y_2))\rVert_{\hk}^2}:=\prkmac,$$
		where $Y_1,Y_2$ are i.i.d.~from $\mu_Y$, and $(Y',\tilde{Y'})$ are defined as follows: Generate $(X',Y',\tilde{Y'})$ as $X'\sim \mu_X$, $Y'|X\sim \mu_{Y|X'}$, $\tilde{Y'}|X'\sim \mu_{Y|X'}$ and $Y'\independent\tilde{Y'}|X'$.\footnote{Note that $Y'$ and $\tilde{Y'}$ are dependent (via $X'$), unlike $Y_1$ and $Y_2$.}
\end{theorem}

\begin{remark}[On condition (S1)]\label{rem:rateassume}
		To better understand~\eqref{eq:artifact} let us first give an alternate expression for $$r(x_1,x_2) =\langle g(x_1),g(x_2)\rangle_{\mathcal{H}_K}, \quad \mbox{where}\quad  g(x):=\E\big[K(R^Y(Y),\cdot) \mid R^X(x)=x\big],$$ i.e., $r(\cdot,\cdot)$ is the inner product between the conditional feature maps (after the rank transformations). Thus,~\eqref{eq:artifact} can be thought of as a Lipschitz assumption on the inner product. A similar condition also appears in~\cite{Azadkia-Chatterjee-2021, DebEtAl-2020}. This condition is needed in our proof technique but we expect~\cref{thm:Consistency} to hold in more generality. 
  
  Let us illustrate how under reasonably practical model assumptions condition (S1) is satisfied. Suppose that one assumes a model where $Y=f(X)+\epsilon$, and $f$ is a Borel measurable function, with the noise variable $\epsilon$ being independent of $X$. Additionally, suppose that $(R^X)^{-1}$ (inverse function of $R^X$), $R^Y$, $f$, and both arguments of the kernel $K(\cdot,\cdot)$ in Theorem~\ref{thm:Consistency} are $\beta$-H\"older continuous.\footnote{The stated assumptions on $R^X$ and $R^Y$ are valid when the distributions of $X$ and $Y$ are supported on closed bounded sets, and their marginal density functions are bounded above and below in their support sets and are $\beta$-H\"older continuous \cite[Corollary~5]{CF19}.} Then we can write
  \begin{align*}
     \Big|r(x_1,x_2) - r(\tilde{x}_1, \tilde{x}_2)\Big| &\leq \mathbb{E}\Big[\big|K\big(R^Y(f((R^X)^{-1}(x_1))+\epsilon_1), f((R^X)^{-1}(x_2))+\epsilon_2)\big)\\ &\,\,\,\,\,\quad - 
     K\big(R^Y(f((R^X)^{-1}(\tilde{x}_1))+\epsilon_1), R^Y(f((R^X)^{-1}(\tilde{x}_2))+\epsilon_2)\big)\big|\Big]\\
      & \leq C(\|x_1-\tilde{x}_1\|^{\beta}+ \|x_2-\tilde{x}_2\|^{\beta})
  \end{align*}
  where $C$ is a constant that depends on the kernel $K$ and the distributions of $X$ and $Y$. The last inequality above follows from the $\beta$-H\"older continuity of $(R^X)^{-1}, R^Y$, $f$ and $K$.
	\end{remark}

\begin{remark}[On condition (S2)]\label{rem:detgraph}
		As the set $\mathbf{\mathcal{H}}_n^{d_1}$, the set of empirical multivariate ranks, i.e., $(\hat{R}_n^X(X_1),\ldots ,\hat{R}_n^X(X_n))$ is some permutation of $\mathbf{\mathcal{H}}_n^{d_1}$, note that $\max_{1\leq i\leq n} d_i$, $\min_{1\leq i\leq n} d_i$, and $\frac{1}{n\min_{i=1}^n d_i}\sum_{e\in \rmgn}|e|^{\beta}$ are all deterministic quantities appearing in condition (S2) above. This should not be confused with the randomness inherent in the construction of $\rgn$. For example, $\mathbf{1}((i,j)\in\rmgn)$, i.e., indicator of the event that $\hat{R}_n^X(X_i)$ and $\hat{R}_n^X(X_j)$ are connected in $\rgn$, is still a random variable.
  \end{remark}

\cref{thm:Consistency} (proved in Section~\ref{sec:pfmain}) shows that $\rkmac$ can be used to construct a test for independence (in addition to being a measure of association) which will be consistent and exactly distribution-free under the null hypothesis of independence.

The above result shows that not only is our measure of association $\rkmac$ distribution-free when $X$ and $Y$ are independent, it converges to a population limit $\prkmac$ as $n$ grows large, under suitable assumptions. We call  $\prkmac$ as the \emph{population measure of association}. Indeed, as the following result shows, $\prkmac$ satisfies many desirable properties that justifies it as a population measure of association between $X$ and $Y$.
\begin{theorem}\label{thm:rankassoc}
If $K(\cdot,\cdot)$ is characteristic and continuous, then $\prkmac$ satisfies the following properties:
\begin{itemize}
		\item[(P1)] $\prkmac  \in [0,1]$.
		\item[(P2)] $\prkmac=0$ if and only if $\mu=\mu_X\otimes \mu_Y$ (i.e., $X$ and $Y$ are independent).
		\item[(P3)] $\prkmac=1$ if and only if $Y=g(X)$, $\mu$ almost everywhere (a.e.), for some measurable function $g:\R^{d_1} \to \R^{d_2}$.
	\end{itemize} 
\end{theorem}

\begin{remark}[Asymmetric nature of $\rkmac$ and $\prkmac$]\label{rem:asymm}
    Note that $\rkmac$, as defined in \eqref{eq:rankkmac}, is not symmetric in the $X_i$s and $Y_j$s. This is by design: The goal of $\rkmac$ is to measure how well \emph{$Y$ can be predicted from $X$}, which can be very different from how well $X$ can be predicted from $Y$. For example, if $Y=f(X)+\epsilon$ for $\epsilon$ independent of $X$, $X\in\R^{d_1}$, $d_1>1$,  and $f:\R^{d_1}\to\R$, then it is conceivable that it is easier to predict $Y$ from $X$ than to predict $X$ from $Y$. Hence it is natural to construct $\rkmac$ in an asymmetric fashion. Of course, it is possible to symmetrize $\rkmac$ by switching the role of $X$, $Y$, and then taking a maximum. It is clear from Theorems \ref{thm:Consistency} and \ref{thm:rankassoc} that this symmetrized version will also converge to a measure of association satisfying (P1), (P2), and a symmetric variant of (P3).
\end{remark}

Next  we show that $\prkmac$ (which arises naturally as the limit of $\rkmac$ in~\cref{thm:Consistency}) is exactly the same as the following population correlation coefficient proposed in~\citet[Equation 1.2]{chatterjee2019new} and~\citet[Equation 2.1]{Azadkia-Chatterjee-2021} where the authors only consider the case $d_2=1$ (also see~\citet{dette2013copula}):
	$$\xi(\mu):=\frac{\int \mathrm{Var}(\mathbb{P}(Y\geq t|X))\,d\mu_Y(t)}{\int \mathrm{Var}(\mathbf{1}(Y\geq t))\,d\mu_Y(t)}$$ for a particular choice of a characteristic kernel. The choice turns out to be the kernel: 
\begin{equation}\label{eq:Dist-Kernel}
 K(y_1,y_2)=|y_1|+|y_2|-|y_1-y_2|.
\end{equation}
Incidentally this is also the kernel used in the construction of distance covariance (see~\cite{Gabor2007}). We state the result formally in the proposition below.
	\begin{prop}\label{prop:Chacon}
		When $d_2=1$ and $K(\cdot,\cdot)$ is the kernel in~\eqref{eq:Dist-Kernel}, then $\prkmac=\xi(\mu)$.
	\end{prop}
The above proposition shows that the family of measures of association proposed in this paper includes that of~\citet{chatterjee2019new} and extends it significantly to general $d_1,d_2\geq 1$ and also a large class of kernel functions. 

\subsection{A rank CLT under null}\label{sec:rankclt}
We present a CLT for $\rkmac$ when $\mu=\mu_X\otimes\mu_Y$. Note that $\rkmac$ is a function of $\hat{R}_n(Y_i)$'s (respectively $\hat{R}_n(X_i)$'s) which are no longer independent among themselves. As a result, proving a CLT for $\rkmac$ is not trivial. Informally speaking, the technique used in this paper is that of a H\'ajek representation (as in~\cite{Shi-Drton-Han-2022}), where we show that the empirical multivariate ranks can be replaced by their population counterparts in $\rkmac$ at a $o_{\mathbb{P}}(1/\sqrt{n})$ cost. As $\rkmac$ does not have a standard $U$-statistic representation, we do not use the explicit form of H\'ajek projections as in~\cite{Shi-Drton-Han-2022}, but opt for a more hands-on method of moments based approach. 

Note that 
\begin{equation}\label{eq:quaden}
\sqrt{n}\,\rkmac=\frac{\nrk}{D_n}
\end{equation}
where $$\;\;D_n:= \frac{1}{n}\sum_{i=1}^n K(\hat{R}_n^Y(Y_i),\hat{R}_n^Y(Y_i))- F_n$$ and 
\begin{equation}\label{eq:recallnrk}
	\nrk:=\sqrt{n}\left(\frac{1}{n}\sum_{i=1}^n\frac{1}{d_i}^{-1}\sum_{j:(i,j)\in\rmgn} K\big(\hat{R}_n^Y(Y_i),\hat{R}_n^Y(Y_j)\big)-F_n\right)
 \end{equation}
with $F_n$ defined in~\eqref{eq:F_n}. Now, $D_n$ does not involve the $X_i$'s and converges to $\E[K(R^Y(Y_1),R^Y(Y_1))]-\E[K(R^Y(Y_1),R^Y(Y_2))]$ by the continuous mapping theorem. Therefore, by Slutsky's theorem, it suffices to establish a pivotal limit distribution for $\nrk$ after a suitable scaling. This is the subject of the following theorem.

\begin{theorem}\label{theo:ranknullclt}
If \eqref{eq:graph1} holds, $K(\cdot,\cdot)$ is continuous, $\mu=\mu_X\otimes\mu_Y$, $\mu_X\in\mathcal{P}_{ac}(\R^{d_1})$ and $\mu_Y\in\mathcal{P}_{ac}(\R^{d_2})$, then $\mathrm{Var}(\nrk)=\mathcal{O}(1)$. Further, with $\theta:=(D,\gamma,\epsilon)\in (0,\infty)^3$ consider the following subclass of graph functionals and measures on $\X \times \Y$ given by:
\begin{small}
    \begin{align*}
		&\mathcal{J}_{\theta} :=\Big\{\tilde{\mathcal{G}}:  \limsup_{n\to\infty}\left(\max_{1\leq i\leq n} \frac{\td_i}{(\log{n})^{\gamma}}+\frac{\max_{i=1}^n \tilde{d}_i}{\min_{i=1}^n \tilde{d}_i}\right)\leq D,\; \mathrm{Var}(\nrk)\geq\epsilon\ \forall n\geq D , \\ & \mathrm{where}  \,
		 (\td_1,\ldots ,\td_n)\; \mathrm{denotes \; the \; degree \; sequence \; of \;}  \trgn:=\tilde{\mathcal{G}}(\hat{R}_n^X(X_1),\ldots ,\hat{R}_n^X(X_n))\Big\}.
		\end{align*}
  \end{small}
		Then, under assumption (S3) from \cref{thm:Consistency}, the following result holds for every fixed $\theta\in (0,\infty)^3$:
		\begin{align}\label{eq:rankmainres}
		\lim\limits_{n\to\infty}\sup_{\tilde{\mathcal{G}}\in \mathcal{J}_{\theta}}\sup_{z\in\R}\Bigg|\P \left(\frac{\nrk}{\hS_n}\leq z\right)-\Phi(z)\Bigg|=0,
		\end{align}
		where 
       $$\tilde{S}_n^2:=\tilde{a}\left(\tilde{g}_1+\tilde{g}_3-\frac{2}{n-1}\right)+\tilde{b}\left(\tilde{g}_2-2\tilde{g}_1-2\tilde{g}_3-\frac{n-5}{n-1}\right)+\tilde{c}\left(\tilde{g}_1+\tilde{g}_3-\tilde{g}_2+\frac{n-3}{n-1}\right)$$ with 
       \begin{align*}
       \tilde{a}&:=\E \left[K^2\big(R^Y(Y_1), R^Y(Y_2)\big)\right],\\ \tilde{b}&:=\E \left[K\big(R^Y(Y_1),R^Y(Y_2)\big) \, K \big(R^Y(Y_1),R^Y(Y_3)\big) \right], \\  \tilde{c}&:=\left(\E \big[K\big(R^Y(Y_1), R^Y(Y_2)\big)\big]\right)^2,
       \end{align*}
       and 
$$\tilde{g}_1:=\frac{1}{n}\sum_{i=1}^n \frac{1}{\tilde{d}_i},\quad \tilde{g}_2:=\frac{1}{n}\sum_{i,j} \frac{T^{\trgn}(i,j)}{\tilde{d}_i \tilde{d}_j}, \quad \tilde{g}_3:=\frac{1}{n}\sum_{(i,j)\in\mathcal{E}(\trgn)} \frac{1}{\tilde{d}_i\tilde{d}_j}.$$
       Here 
       $$T^{\trgn}(i,j):=\sum_{k}\mathbf{1}((i,k)\in\mathcal{E}(\trgn))\, \mathbf{1}((k,j)\in\mathcal{E}(\trgn))$$
       denotes the number of common neighbors of $\hat{R}_n^X(X_i)$ and $\hat{R}_n^X(X_j)$.
\end{theorem}
\begin{remark}[$\tilde S_n^2$ is distribution-free]
    It is worth noting that although $\tilde{S}_n^2$ has a complicated expression, it does not depend on the underlying distributions $\mu_X$ and $\mu_Y$. In particular, $\tilde a, \tilde b, \tilde c$ are fixed constants as $R^Y(Y_1),R^Y(Y_2), $ $R^Y(Y_3)\overset{i.i.d.}{\sim}\mathcal{U}[0,1]^{d_2}$. Moreover, $\tilde g_1, \tilde g_2, \tilde g_3$ are deterministic quantities based on a graph computed from a fixed set of points.
\end{remark}
 
\begin{remark}[Uniformity over $\tilde{\mathcal{G}}$ in Theorem~\ref{theo:ranknullclt}]
One of the strengths of \cref{theo:ranknullclt} is that it is uniform over a large class of graphs (see~\eqref{eq:rankmainres}). This implies that the convergence to normality will hold for data dependent choices of graphs; e.g., when using a $k$-nearest neighbor graph with $k$ chosen as a function of the observed data, provided the choice only grows logarithmically in $n$ with probability converging to $1$ (see the definition of $\mathcal{J}_\theta$).
\end{remark}
 

\section{Discussion and open questions}\label{sec:Disc}
To summarize, in the current manuscript, we have provided a measure of association between two random variables $X\in \X \subset \R^{d_1}$ and $Y\in \Y \subset \R^{d_2}$ which equals $0$ or $1$ depending on whether $Y$ is independent of $X$ or $Y$ is a noiseless measurable function of $X$. We propose a class of tuning-free estimators of this measure which has the added benefit of being exactly distribution-free under the null hypothesis of independence between $X$ and $Y$. Further, we prove a uniform central limit theorem under the null. At a broader conceptual level, our work brings together three different disciplines --- RKHSs, geometric graphs, and OT, in order to achieve the aforementioned set of properties. Our work opens up a number of potential directions to explore in the future. We outline some of them below.

(a) {\bf CLT under alternative}: When $\mu\neq \mu_X\otimes\mu_Y$, the question of proving a limit theorem for $\rkmac$ or $\nrk$ remains open. The main source of difficulty lies in the fact that $\nrk$ no longer has mean $0$ and in fact, the resulting bias decays slower than $n^{-1/2}$. This suggests the need of debiasing $\rkmac$ appropriately. We believe that the works of \cite{berrett2019efficient,berrett2023efficient} could be relevant in this direction.

(b) {\bf Monotonicity}: As $\prkmac=0$ and $1$ capture independence and perfect dependence between $Y$ and $X$, it is natural to expect that the relationship between $Y$ and $X$ should get stronger continuously and monotonically as $\prkmac$ increases from $0$ to $1$. For a variant of $\prkmac$ in the univariate $d_1=d_2=1$ setting, this monotonicity was demonstrated for certain natural examples in \cite[Proposition 4.1]{Auddy-Et-Al-2024}. It would be interesting to establish such monotonicity for general kernels and for $d_1,d_2>1$.

(c) {\bf Local power analysis}: The question of finding the  detection boundary for testing independence using $\nrk$ is of considerable interest. To put it informally, the questions we ask here is: What is the minimum ``distance" between $\mu$ and $\mu_X\otimes \mu_Y$, such that the test for independence using $\rkmac$ has asymptotic power $1$?   Based on the works of \cite{Auddy-Et-Al-2024,Shi-Drton-Han-2022,lin2023boosting}, it is evident that the answer should depend heavily on the choice of the graph functional used in the construction of $\rkmac$, e.g., the choice of $k$ in a $k$-nearest neighbor graph. As in (a) above, the hardness arises again in tightly quantifying the bias in $\nrk$. Here we expect a ``blessing of dimensionality" in the detection boundary beyond $d_2=8$ (see \citet{BhaswarAoS19} for a related result).    

(d) {\bf Numerical performance}: We do not pursue a thorough simulation study comparing the power of the distribution-free test of independence based on $\rkmac$ with other competing methods here. However, based on recent experimental observations made for other graph-based and rank-based tests, there are a number of interesting practical questions one may ask. For instance, graph-based tests usually outperform competitors when the relationship between $Y$ and $X$ is sinusoidal (see \cite{berrett2017nonparametric,chatterjee2019new}) while multivariate rank-based tests are known for their robustness against outliers and heavy-tailed data (see~\cite{Deb19,shi2020rate}). As $\rkmac$ uses both, one may ask: Does it inherit the best of both worlds? Also, in \cite[Section D.5]{deb2021pitman}, the authors show that certain multivariate rank-based tests have competitive performance in high-dimensional problems as well. It would be interesting to explore if $\rkmac$ exhibits similar behavior.

(e) {\bf Extension beyond uniform reference distribution}: Throughout this paper we have assumed that the reference distributions needed for defining the multivariate ranks are uniform on the $d_1$ and $d_2$-dimensional unit hypercube. This choice is made mostly for simplicity and from a methodological perspective any absolutely continuous reference distribution would have sufficed (see e.g., \cite{deb2021pitman,shi2021center,hallin2023efficient} for some alternate choices of the reference distribution). Our proofs would work verbatim for any other compactly supported reference distribution. However, for reference distributions with unbounded supports (e.g., multivariate normal) some truncation arguments would be needed. The benefits of using other reference distributions have a rich literature dating back to \cite{hodges1956efficiency,chernoff1958asymptotic} and more recently in \cite{deb2021pitman,shi2021center}. It remains open to explore whether such gains can be had by choosing other reference distributions for constructing $\rkmac$.

\section{Proofs of our main results}\label{sec:pfmain}

	\subsection{Proof of Theorem~\ref{thm:Consistency}}\label{pf:Consistency}
	
	First we state an important result about the consistency of the empirical multivariate ranks  (as in~\cref{def:empquanrank}), i.e., $\hat{R}_n(\cdot)$ yields a consistent estimate of $R(\cdot)$ (see~\cite{Deb19} for a proof).
	\begin{prop}\label{prop:consisrank}
		Suppose $\mu\in\mathcal{P}_{ac}(\R^d)$ and $n^{-1}\sum_{i=1}^n \delta_{h_i^d}\overset{w}{\longrightarrow}\nu$. Then $n^{-1}\sum_{i=1}^n \lVert \hat{R}_n(X_i)-R(X_i)\rVert_2\overset{a.s.}{\longrightarrow}0$.
	\end{prop}

	\emph{Proof of (i)}: By~\cite[Proposition 2.2 (ii)]{Deb19}, $(\hat{R}_n^X(X_1),\ldots ,\hat{R}_n^X(X_n))$ (or $(\hat{R}_n^Y(Y_1),\ldots ,\hat{R}_n^Y(Y_n))$) is distributed uniformly over all $n!$ permutations of the set $\mathbf{\mathcal{H}}_n^{d_1}$ (or $\mathbf{\mathcal{H}}_n^{d_2}$). When $\mu=\mu_X\otimes\mu_Y$, clearly $(\hat{R}_n^X(X_1),\ldots ,\hat{R}_n^X(X_n))$ and $(\hat{R}_n^Y(Y_1),\ldots ,\hat{R}_n^Y(Y_n))$ are independent and the joint distribution is pivotal. As $\rkmac$ is a function of $(\hat{R}_n^X(X_1),\ldots ,\hat{R}_n^X(X_n),\hat{R}_n^Y(Y_1),\ldots ,\hat{R}_n^Y(Y_n))$, it also has a pivotal distribution.
	
	\emph{Proof of (ii)}: As $n^{-1}\sum_{i=1}^n \delta_{h_i^{d_2}}$ converges weakly to $\mathcal{U}[0,1]^{d_2}$, the following conclusions are easy consequences of the Portmanteau Theorem:
	\begin{align}\label{eq:Promitexp}
	&\frac{1}{n}\sum_{i=1}^n K(\hat{R}_n^Y(Y_i),\hat{R}_n^Y(Y_i))\overset{n\to\infty}{\longrightarrow}\E K(R^Y(Y_1),R^Y(Y_1)),\nonumber  \\&
	\frac{1}{n(n-1)}\sum_{i\neq j} K(\hat{R}_n^Y(Y_i),\hat{R}_n^Y(Y_j))\overset{n\to\infty}{\longrightarrow}\E K(R^Y(Y_1),R^Y(Y_2)).
	\end{align}
	In the above displays, both terms on the right hand side are deterministic because both of them are permutation-invariant and $(\hat{R}_n(Y_1),\ldots ,\hat{R}_n(Y_n))$ is some permutation of the fixed set $\mathcal{H}_n^{d_2}$. Consequently the above convergence is a deterministic result.
	By an application of~\cref{prop:consisrank} and \eqref{eq:graph1}, it further suffices to show that:
	$$\mathcal{Z}_n:=\frac{1}{n}\sum_{i}d^{-1}_i\sum_{j:(i,j)\in \rmgn}K(R^Y(Y_i), R^Y(Y_j))\overset{\mathbb{P}}{\longrightarrow}\E K(R^Y(Y'),R^Y(\tilde{Y'})).$$
	
	In order to establish the above, by an application of Chebyshev's inequality, it suffices to show that:
	\begin{enumerate}
		\item[(a)] $\E \mathcal{Z}_n^2\overset{n\to\infty}{\longrightarrow} \left[\E K(R^Y(Y'),R^Y(\tilde{Y'}))\right]^2.$
		\item[(b)] $\E \mathcal{Z}_n\overset{n\to\infty}{\longrightarrow} \E K(R^Y(Y'),R^Y(\tilde{Y'})).$
	\end{enumerate}
	
	We first prove (a). For the sake of simplicity, we write 
	\begin{align*}
	\mathbb{E}\mathcal{Z}_n^2 & = (\mathbf{I})+ (\mathbf{II}) + (\mathbf{III})
	\end{align*}
	where 
	\begin{align*}
	(\mathbf{I}) &:= \frac{1}{n^2}\E\Bigg[\sum_{i}\sum_{j:(i,j)\in \rmgn}\frac{1}{d_i}\left(\frac{1}{d_i}+\frac{1}{d_j}\right)\mathbb{E}\Big[ K^2(R^Y(Y_i), R^Y(Y_j)) \mid R^X(X_i),R^X(X_j)\Big]\Bigg] \\
	(\mathbf{II}) &:= \frac{1}{n^2}\E\Bigg[\sum_{i}\sum_{j\neq k: (i,j),(i,k)\in \rmgn}\left(\frac{1}{d_i}+\frac{1}{d_j}\right)\left(\frac{1}{d_i}+\frac{1}{d_k}\right)\\& \quad \times\mathbb{E}\Big[K(R^Y(Y_i), R^Y(Y_j))\,K(R^Y(Y_i), R^Y(Y_k)) \mid R^X(X_i),R^X(X_j),R^X(X_k)\Big]\Bigg]\\
	(\mathbf{III}) &:=\frac{1}{n^2} \E\Bigg[\sum_{i\neq j} \sum_{k\neq \ell: (k,i), (\ell,j)\in\rmgn}\frac{1}{d_id_j} \mathbb{E}\Big[K(R^Y(Y_i),R^Y(Y_k))\\ &\qquad \qquad \times K(R^Y(Y_j), R^Y(Y_\ell)) \mid  R^X(X_i),R^X(X_j),R^X(X_k),R^X(X_l)\Big]\Bigg].  
	\end{align*}
	We now claim and prove that 
	\begin{align}\label{eq:ppconv}
	(\mathbf{III}) \to \Big(\mathbb{E}\Big(\mathbb{E}\Big[K\big(R^Y(\tilde{Y}_1),R^Y(\tilde{Y}_2)\big) \mid X\Big]\Big)\Big)^2
	\end{align}
	as $n\to\infty$, where $\tilde{Y}_1$ and $\tilde{Y}_2$ are independent samples from the conditional distribution $Y\mid X$. 
	For this, we define 
	\begin{align*}
	(\widetilde{\mathbf{III}}) := \frac{1}{n^2}\sum_{i\neq j} \sum_{k\neq \ell: (k,i),(\ell,j)\in\rmgn} &\frac{1}{d_i}\frac{1}{d_j}\mathbb{E}\Big[\mathbb{E}\Big[K\big(R^Y(Y_i),R^Y(Y^{\prime}_i)\big)\mid R^X(X_i)\Big]\\&\times \mathbb{E}\Big[K\big(R^Y(Y_j),R^Y(Y^{\prime}_j)\big)\mid R^X(X_j)\Big]\Big]
	\end{align*}
	where $Y^{\prime}_i$ (resp. $Y^{\prime}_j$) is a sample from the conditional distribution $Y\mid X_i$ (resp. $Y\mid X_j$) independent of $Y_i$ and $Y_j$.  
	A simple application of the fact that $K(\cdot,\cdot)$ is continuous a.e., $R^Y(\cdot)$ is bounded coordinate-wise, shows that:
	\begin{align}
	&\;\;\;|(\mathbf{III}) -(\widetilde{\mathbf{III}})|\nonumber \\ &\lesssim C\frac{\max_{j=1}^n d_j}{n\min_{j=1}^n d_j}\sum_{i}\frac{1}{d_i}\mathbb{E}\Big[\sum_{k:(k,i)\in \rmgn}\big|  \mathbb{E}\big[K(R^Y(Y_i),R^Y(Y_k))\mid R^X(X_i),R^X(X_k)\big]\nonumber \\ &\qquad  - \mathbb{E}\big[K(R^Y(Y_i),R^Y(Y^{\prime}_i))\mid R^X(X_i)\big]\big|\Big].  \label{eq:Error}
	\end{align}
	We show that the right hand side of the above inequality converges to $0$ as $n \to \infty$. By the assumption \eqref{eq:artifact},  $\mathbb{E}\big[K\big(R^Y(Y_1), R^Y(Y_2)\big)\mid R^X(X_1)=x_1,R^X(X_2)=x_2\big]$ is uniformly $\beta$-H{\"o}lder continuous w.r.t. $x_1,x_2$. Applying this to \eqref{eq:Error} shows 
	\begin{align}
	&\;\;\;\;\limsup_{n\to \infty}\text{r.h.s. of \eqref{eq:Error}}\nonumber \\ &\lesssim \limsup_{n\to \infty}\frac{\max_{j=1}^n d_j}{n\min_{j=1}^n d_j}\sum_{i=1}^{n}\frac{1}{d_i}\mathbb{E}\Big[\sum_{k:(k,i)\in \rmgn} \big\|R^X(X_i)- R^X(X_k)\big\|^{\beta} \Big]\nonumber \\
	&\leq \limsup_{n\to \infty} \frac{\max_{j=1}^n d_j}{n\min_{j=1}^n d_j} \mathbb{E}\Big[\sum_{i=1}\frac{1}{d_i}\sum_{k:(k,i)\in\rmgn}\big\|\hat{R}_n^X(X_i)-\hat{R}_n^X(X_k)\big\|^{\beta}\Big] . \label{eq:errorbd}
	\end{align} 
	where the last inequality follows once again from~\cref{prop:consisrank}. Notice that $\sum_{i=1}\sum_{k:(k,i)\in\rmgn}\big\|\hat{R}_n^X(X_i)-\hat{R}_n^X(X_k)\big\|^{\beta}=\sum_{e\in\rmgn} |e|^{\beta}$ where $|e|$ denotes the length of the edge $e$. The right hand side of the above display converge to $0$ by \eqref{eq:graph2}.   
	This shows $|(\mathbf{III})-(\widetilde{\mathbf{III}})|$ converges to $0$ as $n$ tends to $\infty$. On other hand, we define 
	\begin{align*}
	(\widetilde{\mathbf{I}}) &:= \frac{1}{n^2} \E\Bigg[\sum_{i=1}^{n} \sum_{j:(j,i)\in \rmgn} \Bigg(\frac{1}{d^2_i}\Big[\big(K(R^Y(Y_i), R^Y(Y^{\prime}_i))\big)^2\Big]\\& \qquad \qquad + \frac{1}{d_id_j}\mathbb{E}\Big[K(R^Y(Y_i), R^Y(Y^{\prime}_i))\, K(R^Y(Y_j), R^Y(Y^{\prime}_j))\Big]\Bigg)\Bigg],\\
	\end{align*}
and
 \begin{align*}
&(\widetilde{\mathbf{II}}) :=\frac{1}{n^2}\E\Bigg[\sum_{i=1}^{n} \sum_{(j,i),(k,i)\in\rmgn}\Bigg[ \frac{1}{d_jd_k}K(R^Y(Y_k), R^Y(Y^{\prime}_k))\, K(R^Y(Y_j), R^Y(Y^{\prime}_j))\\ & \,\, \quad +\frac{1}{d_i^2}K^2(R^Y(Y_i),R^Y(Y_i'))+\frac{2}{d_id_j}K(R^Y(Y_i),R^Y(Y_i'))\, K(R^Y(Y_j),R^Y(Y_j'))\Bigg]\Bigg].
\end{align*}
	Note that 
	\begin{align}
	(\widetilde{\mathbf{I}}) +(\widetilde{\mathbf{II}}) + (\widetilde{\mathbf{III}}) &= \frac{1}{n^2} \mathbb{E}\Big[\Big(\sum_{i=1}^n K(R^Y(Y_i), R^Y(Y^{\prime}_i))\Big)^2\Big]\nonumber \\& \to \Big(\mathbb{E}\big[K(R^Y(Y_1), R^Y(Y^{\prime}_1))\big]\Big)^2 \label{eq:sumconv}
	\end{align}
	where the last convergence follows from the strong law of large numbers and the dominated convergence theorem. Owing to the fact that \eqref{eq:graph1} holds and $K(\cdot,\cdot)$ is continuous a.e., it is easy to check that $(\mathbf{I}),(\mathbf{II}),(\widetilde{\mathbf{I}}),(\widetilde{\mathbf{II}})$ are converging to $0$ as $n \to \infty$. Combining this with \eqref{eq:ppconv} and \eqref{eq:sumconv} completes the proof of (a). 
	
	Now we move on to (b). By the towering property of the conditional expectation, we have 
	\begin{align}
	\E \mathcal{Z}_n&=\mathbb{E}\Big[\frac{1}{n}\sum_{i}d^{-1}_i\sum_{j:(i,j)\in \mathcal{E}(\mathcal{G}_n)}K(R^Y(Y_i), R^Y(Y_j))\Big]\nonumber \\
	& = \mathbb{E}\Big[\frac{1}{n}\sum_{i}d^{-1}_i\sum_{j:(i,j)\in \mathcal{E}(\mathcal{G}_n)}\mathbb{E}\big[K(R^Y(Y_i), R^Y(Y_j))\mid R^X(X_i),R^X(X_j)\big]\Big].\label{eq:Tower}
	\end{align}
	By the $\beta$-H{\"o}lder continuity of $\mathbb{E}\big[K(R^Y(Y_i), R^Y(Y_j))\mid R^X(X_i)=x, R^X(X_j)=y\big]$ as a function of $x$ and $y$ (see \eqref{eq:artifact}), there exists $C>0$ such that 
	\begin{align*}
	\Big|\,&\mathbb{E}\big[K(R^Y(Y_i), R^Y(Y_j))\mid R^X(X_i), R^X(X_j)\big]\\& -\mathbb{E}\big[K(R^Y(Y_i), R^Y(Y^{\prime}_j))\mid R^X(X_i)\big] \, \Big| \lesssim \|R^X(X_i)- R^X(X_j)\|^{\beta}
	\end{align*}
	where $Y_i,Y^{\prime}_i$ are two independent samples from the conditional distribution of $Y$ given $X_i$. Applying the above inequality, we get 
	\begin{align*}
	\Big|&\text{r.h.s. of \eqref{eq:Tower}}- \frac{1}{n}\mathbb{E}\big[\sum_{i=1}^{n}K(R^Y(Y_i), R^Y(Y^{\prime}_i))\big]\Big|\\& \lesssim \mathbb{E}\Big[\frac{1}{n}\sum_{i=1}^{n}\sum_{j:(i,j)\in \mathcal{E}(\mathcal{G}_n)}\frac{1}{d_i}\|R^X(X_i)- R^X(X_j)\|^{\beta}\Big].
	\end{align*} 
	The proof of (b) can now be completed using the same steps as those in the proof of (a) starting from~\eqref{eq:errorbd}.    	
	\qed 

 \subsection{Proof of Theorem \ref{thm:rankassoc}}\label{pf:rankassoc}

\emph{Proof of (P1)}: This follows from \cite[Theorem 2.1]{DebEtAl-2020}.

\noindent \emph{Proof of (P2)}: By \cref{prop:Mccan}, there exists $Q^Y(\cdot)$ such that $Q^Y(R^Y(y))=y$ a.e. $\mu_Y$. By \cite[Theorem 2.1]{DebEtAl-2020}, it follows that $\prkmac=0$ if and only if $R^{Y}(Y)$ and $X$ are independent. As $Q^Y(R^Y(Y))=Y$, we further have that $R^{Y}(Y)$ and $X$ are independent if and only if $Y$ and $X$ are independent. The conclusion then follows. 

\noindent \emph{Proof of (P3)}: By \cite[Theorem 2.1]{DebEtAl-2020}, it follows that $\prkmac=1$ if and only if there exists a measurable function $h(\cot)$ such that $R^{Y}(Y)=h(X)$ a.e. $\mu$. Now $R^Y(Y)=h(X)\Leftrightarrow Y=Q^Y(h(X))$. The conclusion then follows by choosing $g\equiv Q^Y\circ h$.\qed

\subsection{Proof of Proposition~\ref{prop:Chacon}}\label{pf:Chacon}
	
	The crucial observation in this proof is that $K(y_1,y_2)=|y_1|+|y_2|-|y_1-y_2|=\min\{y_1,y_2\}$ for $y_1,y_2\in [0,\infty)$. Now, when $d_2=1$, $R^Y(\cdot)$ is simply the cumulative distribution function of $Y$; we will call it $F_Y$ to stick with conventional notation. Next note that, $\rpk$ can be simplified as:
	\begin{equation}\label{eq:temprpk}
	\rpk=\frac{\E\min\{F_Y(Y'),F_Y(\tilde{Y'})\}-\E\min\{F_Y(Y_1),F_Y(Y_2)\}}{1/2-\E\min\{F_Y(Y_1),F_Y(Y_2)\}}.
	\end{equation}
	As $\min\{a,b\}=\int_0^{\infty} \mathbf{1}(t\leq a)\mathbf{1}(t\leq b)\,dt$ for $a,b\in\R$, an application of the dominated convergence theorem yields: $$\E\min\{F_Y(Y'),F_Y(\tilde{Y'})\}=\int \E [\left(\P(F_Y(Y)\geq t|X)\right)^2]\,dt=\int \E[\left(\P(Y\geq t|X)\right)^2]\,d\mu_Y(t)$$ and similarly, $$\E\min\{F_Y(Y_1),F_Y(Y_2)\}=\int \left(\P(Y\geq t)\right)^2\,d\mu_Y(t).$$ Plugging the above expressions into the expression of $\rpk$ in~\eqref{eq:temprpk} gives: $$\rpk =\frac{\int\left(\E[\left(\P(Y\geq t|X)\right)^2]-\left(\P(Y\geq t)\right)^2\right)\,d\mu_Y(t)}{\int\left(\P(Y\geq t)-\left(\P(Y\geq t)\right)^2\right)\,d\mu_Y(t)}=\xi(\mu).$$ This completes the proof.
	\qed
	
	\subsection{Proof of Theorem \ref{theo:ranknullclt}}\label{pf:ranknullclt}
Note the definitions of $\tilde{a}$, $\tilde{b}$, $\tilde{c}$, $\tilde{g}_1$, $\tilde{g}_2$, $\tilde{g}_3$ from the statement of \cref{theo:ranknullclt}. 
For this proof, we will also borrow some elements used in the proof of \cite[Theorem 4.1]{DebEtAl-2020}. 
Define $C_n:=(n(n-1))^{-1}\sum_{i\neq j} K(\hat{R}_n^Y(Y_i),\hat{R}_n^Y(Y_j))$ and recall that $C_n$ is a deterministic quantity as was explained in the comment after~\eqref{eq:Promitexp}. Also $\E K(\hat{R}_n^Y(Y_1),\hat{R}_n^Y(Y_2))=C_n$. Let $\mathcal{F}_n:=\sigma(X_1,\ldots ,X_n)$, i.e., the $\sigma$-field generated by $(X_1,\ldots ,X_n)$. Consequently, note that:
	\begin{align}\label{eq:rcenter}
	\E[\nrk|\mathcal{F}_n]&=\sqrt{n}\left(\frac{1}{n}\sum_{i=1}^n \frac{1}{\tilde{d}_i}\sum_{j:(j,i)\in\rmgn}\E [K(\hat{R}_n^Y(Y_i),\hat{R}_n^Y(Y_j))]-C_n\right)\nonumber \\&=\sqrt{n}C_n\left(\frac{1}{n}\sum_{i=1}^n \frac{1}{\tilde{d}_i}\sum_{j:(j,i)\in\rmgn} 1-1\right)=0.
	\end{align}
	By the same calculation as in \cite[Equation (C.12)]{DebEtAl-2020}, we get:
	\begin{align}\label{eq:rankmainvareq}
	\mbox{Var}(\nrk|\mathcal{F}_n)=\left(\tilde{g}_1+\tilde{g}_3-\frac{2}{n-1}\right)(\hat{a}-2\hat{b}+\hat{c})+(\tilde{g}_2-1)(\hat{b}-\hat{c}),
	\end{align}
	where 
	\begin{align*}
	&\hat{a}:=\frac{1}{n(n-1)}\sum_{(i,j)\ \mathrm{distinct}} K^2(\hat{R}_n^Y(Y_i),\hat{R}_n^Y(Y_j))\\ &\hat{b}:=\frac{1}{n(n-1)(n-2)}\sum_{(i,j,l)\ \mathrm{distinct}}K(\hat{R}_n^Y(Y_i),\hat{R}_n^Y(Y_j))K(\hat{R}_n^Y(Y_i),\hat{R}_n^Y(Y_l))\\ &\hat{c}:=\frac{1}{n(n-1)(n-2)(n-3)}\sum_{(i,j,l,m)\ \mathrm{distinct}} K(\hat{R}_n^Y(Y_i),\hat{R}_n^Y(Y_j))K(\hat{R}_n^Y(Y_l),\hat{R}_n^Y(Y_m)).
	\end{align*}
	Now $\hat{a},\hat{b},\hat{c}$ are clearly $\mathcal{O}(1)$ and $\tilde{g}_1,\tilde{g}_2,\tilde{g}_3$ are $\mathcal{O}(1)$ by using ~\cite[Equation (C.13)]{DebEtAl-2020}. Next observe that the right hand side of~\eqref{eq:rankmainvareq} is a deterministic quantity. Therefore by combining~\eqref{eq:rcenter}~and~\eqref{eq:rankmainvareq}, we get $\mbox{Var}(\nrk)=\mathcal{O}(1)$.
	
	In order to establish the CLT, define
	\begin{equation}\label{eq:hajekpop}
	\nrp:=\sqrt{n}\left(\frac{1}{n}\sum_{i=1}^n \frac{1}{\tilde{d}_i}\sum_{\substack{j:(j,i)\\ \in\rmgn}} K(R^Y(Y_i),R^Y(Y_j))-\frac{\sum_{i\neq j} K(R^Y(Y_i),R^Y(Y_j))}{n(n-1)}\right).
	\end{equation}
	
 It then suffices to show the following:
 \begin{align}\label{eq:toshow1}
\sup_{\tilde{\mathcal{G}}\in\mathcal{J}_{\theta}}\Bigg|\frac{\mbox{Var}(\nrp)}{\mbox{Var}(\nrk)}-1\Bigg|\overset{n\to\infty}{\longrightarrow}0,
\end{align}
\begin{align}\label{eq:toshow2} \sup_{\tilde{\mathcal{G}}\in\mathcal{J}_{\theta}}\E\left[\left(\nrp-\nrk\right)^2\right]\overset{n\to\infty}{\longrightarrow}0,
 \end{align}
\begin{align}\label{eq:toshow3}
\sup_{\tilde{\mathcal{G}}\in \mathcal{J}_{\theta}}\sup_{x\in\R} \bigg|\P\left(\frac{\nrp}{\sqrt{\mbox{Var}(\nrp)}}\le x\right)-\Phi(x)\bigg|\overset{n\to\infty}{\longrightarrow} 0.
\end{align}

\emph{Proof of \eqref{eq:toshow1}.}	As in \cite[Equation (C.12)]{DebEtAl-2020}, we have:
	\begin{equation}\label{eq:popvar}
	\mbox{Var}(\nrp)=\E(\nrp)^2=\left(\tilde{g}_1+\tilde{g}_3-\frac{2}{n-1}\right)(\tilde{a}-2\tilde{b}+\tilde{c})+(\tilde{g}_2-1)(\tilde{b}-\tilde{c}).
	\end{equation}
	Therefore, 
	\begin{align*}
	&\sup_{\tilde{\mathcal{G}}\in\mathcal{J}_{\theta}}\Bigg|\frac{\mbox{Var}(\nrp)}{\mbox{Var}(\nrk)}-1\Bigg|\lesssim \sup_{\tilde{\mathcal{G}}\in\mathcal{J}_{\theta}}\Bigg[\left(\tilde{g}_1+\tilde{g}_3-\frac{2}{n-1}\right)(\tilde{a}-\hat{a}-2\tilde{b}+2\hat{b}+\tilde{c}-\hat{c})\\ &\qquad \qquad +(\tilde{g}_2-1)(\tilde{b}-\hat{b}-\tilde{c}+\hat{c})\Bigg]\overset{n\to\infty}{\longrightarrow}0.
	\end{align*}
 The conclusion above follows by noting that $\hat{a}\to \tilde{a}$, $\hat{b}\to\tilde{b}$, and $\hat{c}\to \tilde{c}$, as $n\to\infty$, by assumption (S3). 

 \vspace{0.1in}
 
 \emph{Proof of \eqref{eq:toshow2}.} Let us first introduce the notion of the \emph{resampling distribution}. Note that when $\tilde{\mu}=\tilde{\mu}_X\otimes\tilde{\mu}_Y$, the joint distribution of $(X_1,Y_1),\ldots ,(X_n,Y_n)$ is the same as the joint distribution of $(X_1,Y_{\sigma(1)}),\ldots ,(X_n,Y_{\sigma(n)})$ where $\sigma$ is a random permutation of the set $\{1,2,\ldots,n\}$ which is drawn independent of the $X_i$'s and $Y_i$'s. 
		
		The expressions for $\E(\nrk)^2$ and $\E(\nrp)^2$ has already been presented in~\eqref{eq:rankmainvareq}~and~\eqref{eq:popvar} respectively. Therefore, in the sequel, we will only focus on the term $\E[\nrk\nrp]$ where we will use the \emph{resampling distribution} as discussed above. Towards this direction, let us define:
		\begin{align*}
		&\underline{a}:=\frac{1}{n(n-1)}\sum_{(i,j)\ \mathrm{distinct}} K(\hat{R}_n^Y(Y_i),\hat{R}_n^Y(Y_j))K(R^Y(Y_i),R^Y(Y_j))\\ &\underline{b}:=\frac{1}{n(n-1)(n-2)}\sum_{(i,j,l)\ \mathrm{distinct}}K(\hat{R}_n^Y(Y_i),\hat{R}_n^Y(Y_j))K(R^Y(Y_i),R^Y(Y_l))\\ &\underline{c}:=\frac{1}{n(n-1)(n-2)(n-3)}\sum_{(i,j,l,m)\ \mathrm{distinct}} K(\hat{R}_n^Y(Y_i),\hat{R}_n^Y(Y_j))K(R^Y(Y_l),R^Y(Y_m)).
		\end{align*}
		Also to simplify notation, we will use the symbol $j\sim i$ for $(j,i)\in\rmgn$ and let $\hat{Z}_i:=\hat{R}_n^Y(Y_i)$ and $Z_i:=R^Y(Y_i)$. With the above notation, observe that $\E[\nrk\nrp]$ can be simplified as follows: 
		\begin{align*}
		&\E\Bigg[\left(\frac{1}{n}\sum_{i=1}^n \frac{1}{\tilde{d}_i}\sum_{j\sim i} K(\hat{R}_n^Y(Y_{\sigma(i)}),\hat{R}_n^Y(Y_{\sigma(j)}))\right)\left(\frac{1}{n}\sum_{i=1}^n \frac{1}{\tilde{d}_i}\sum_{j\sim i} K(R^Y(Y_{\sigma(i)}),R^Y(Y_{\sigma(j)}))\right)\\ &-\left(\frac{1}{n}\sum_{i=1}^n \frac{1}{\tilde{d}_i}\sum_{j\sim i} K(\hat{R}_n^Y(Y_{\sigma(i)}),\hat{R}_n^Y(Y_{\sigma(j)}))\right)\left(\frac{1}{n(n-1)}\sum_{i\neq j} K(R^Y(Y_i),R^Y(Y_j))\right)\\ &-\left(\frac{1}{n}\sum_{i=1}^n \frac{1}{\tilde{d}_i}\sum_{j\sim i} K(R^Y(Y_{\sigma(i)}),R^Y(Y_{\sigma(j)}))\right)\left(\frac{1}{n(n-1)}\sum_{i\neq j} K(\hat{R}_n^Y(Y_i),\hat{R}_n^Y(Y_j))\right)\\ &+\left(\frac{1}{n(n-1)}\sum_{i\neq j}K(\hat{R}_n^Y(Y_{\sigma(i)}),\hat{R}_n^Y(Y_{\sigma(j)}))\right)\left(\frac{1}{n(n-1)}\sum_{i\neq j} K(R^Y(Y_{\sigma(i)}),R^Y(Y_{\sigma(j)}))\right)\Bigg]\\ &\overset{(*)}{=}(\tilde{g}_1+\tilde{g}_3)\E\left[K(\hat{Z}_{\sigma(1)},\hat{Z}_{\sigma(2)})K(Z_{\sigma(1)},Z_{\sigma(2)})\right]+(3+\tilde{g}_2-2\tilde{g}_1-2\tilde{g}_3)\E\Big[K(\hat{Z}_{\sigma(1)},\hat{Z}_{\sigma(2)})\\ &\quad K(Z_{\sigma(1)},Z_{\sigma(3)})\Big]+(n-3+\tilde{g}_1+\tilde{g}_2+\tilde{g}_3)\E\left[K(\hat{Z}_{\sigma(1)},\hat{Z}_{\sigma(2)})K(Z_{\sigma(3)},Z_{\sigma(4)})\right]\\ &\quad-\frac{2}{n-1}\cdot\E\left[K(\hat{Z}_{\sigma(1)},\hat{Z}_{\sigma(2)})K(Z_{\sigma(1)},Z_{\sigma(2)})\right]-\frac{4(n-2)}{n-1}\cdot\E\left[K(\hat{Z}_{\sigma(1)},\hat{Z}_{\sigma(2)})K(Z_{\sigma(1)},Z_{\sigma(3)})\right]\\ &\quad-\frac{(n-2)(n-3)}{n-1}\cdot\E\left[K(\hat{Z}_{\sigma(1)},\hat{Z}_{\sigma(2)})K(Z_{\sigma(3)},Z_{\sigma(4)})\right]\\ &=\left(\tilde{g}_1+\tilde{g}_3-\frac{2}{n-1}\right)(\underline{a}-2\underline{b}+\underline{c})+(\tilde{g}_2-1)(\underline{b}-\underline{c}).
		\end{align*}
		Here $(*)$ follows from \cite[Equations (C.9), (C.10), and (C.11)]{DebEtAl-2020}.
		By Assumption (S3), $\hat{a},\underline{a}$ converge to $\tilde{a}$ in $L^1$; same holds for $\hat{b},\underline{b},\tilde{b}$ and $\hat{c},\underline{c},\tilde{c}$. Also $\hat{a},\underline{a},\tilde{a},\hat{b},\underline{b},\tilde{b},\hat{c},\underline{c},\tilde{c}$ do not depend on the graph functional $\mathcal{G}$. Therefore, the following holds:
		\begin{align*}
		&\sup_{\tilde{\mathcal{G}}\in\mathcal{J}_{\theta}}\E\left[\left(\nrk-\nrp\right)\right]^2=\sup_{\tilde{\mathcal{G}}\in\mathcal{J}_{\theta}}\Bigg[\left(\tilde{g}_1+\tilde{g}_3-\frac{2}{n-1}\right)\{\hat{a}-2\underline{a}+\tilde{a}+4\underline{b}-2\tilde{b}\\ &\qquad \qquad  -2\hat{b}+\hat{c}-2\underline{c}+\tilde{c} \} +(\tilde{g}_2-1)(\hat{b}-2\underline{b}+\tilde{b}-\hat{c}+2\underline{c}-\tilde{c})\Bigg]\overset{n\to\infty}{\longrightarrow}0.
		\end{align*}
	This completes the proof of \eqref{eq:toshow2}.

 \vspace{0.1in}

 \emph{Proof of \eqref{eq:toshow3}.} The conclusion follows directly from \cite[Theorem 4.1]{DebEtAl-2020}.
 \qed

\bibliographystyle{chicago}
\bibliography{OT,template,References}

\end{document}